\documentclass[12pt]{amsart}
\usepackage{amssymb}

\newtheorem{Thm}{Theorem}[section]

\newtheorem{Cor}[Thm]{Corollary}
\newtheorem{Lem}[Thm]{Lemma}
\newtheorem{Prop}[Thm]{Proposition}

\theoremstyle{definition}

\theoremstyle{remark}

\def\ldots{\mathinner{\ldotp\ldotp\ldotp}}
\def\ldots{\mathinner{\cdotp\cdotp\cdotp}}

\def \cal{\mathcal}
\def \Bbb{\mathbb}

\begin{document}

\title{On subspaces of $c_0$ and extension of operators into
$C(K)$-spaces
}
\author{N. J. Kalton}
\address{Department of Mathematics \\
University of Missouri-Columbia \\
Columbia, MO 65211
}

\email{nigel@math.missouri.edu}


\thanks{The author was supported by NSF grant DMS-9870027}

\begin{abstract} Johnson and Zippin recently showed that if $X$ is a
weak$^*$-closed subspace of $\ell_1$ and $T:X\to C(K)$ is any bounded
operator then $T$ can extended to a bounded operator $\tilde T:\ell_1\to
C(K).$  We give a converse result: if $X$ is a subspace of $\ell_1$ so
that $\ell_1/X$ has a (UFDD) and every operator $T:X\to C(K)$ can be
extended to $\ell_1$ then there is an automorphism $\tau$ of $\ell_1$ so
that $\tau(X)$ is weak$^*$-closed.  This result is proved by
studying subspaces of $c_0$ and several different
characterizations of such subspaces are given.

\end{abstract}

\maketitle

\section{Introduction} In \cite{JZ}, Johnson and Zippin proved an
extension theorem for operators into $C(K)$-spaces:

\begin{Thm}\label{JZth} Let $X$ be a weak$^*$-closed subspace of $\ell_1$
(considered as the dual of $c_0$) and let $T:X\to C(K)$ be a bounded
operator.  Then $T$ has  an extension $\tilde T:\ell_1\to
C(K).$\end{Thm}

Note that this implies the same conclusion for any subspace $X$ so that
$\ell_1/X$ is isomorphic to the dual of a subspace of $c_0$ (using
results of \cite{LR}).  The aim of this paper is to prove a partial
converse
result to the Johnson-Zippin theorem.  We show that if $X$ is a subspace
of $\ell_1$ such that every bounded operator $T:X\to C[0,1]$ can be
extended and if, additionally, $\ell_1/X$ has a (UFDD) (unconditional
finite-dimensional decomposition) then
$\ell_1/X$
is isomorphic to the dual of a subspace of $c_0$, and hence there is an
automorphism $\tau$ of $\ell_1$ such that $\tau(X)$ is weak$^*$-closed.
The hypothesis on $X$ can be weakened a little: it suffices that
$\ell_1/X$ be the dual of space which embeds in a space with a (UFDD).

 The technique of proof depends heavily on ideas developed in \cite{GKL},
where subspaces of $c_0$ are characterized in terms of properties of
norms.  We also use ideas from \cite{GKL2} where trees are used to obtain
renormings, to obtain a characterization of subspaces of $c_0$ in terms
of properties of trees in the dual.

If $X$ is a subspace of $\ell_1$ which satisfies the conclusion of
Theorem \ref{JZth} we show that $\ell_1/X$ has a property we call the
very strong Schur property (the strong Schur property was considered
first for subspaces of $L_1$ by Rosenthal \cite{R}; see also \cite{BR}).
In the presence of some unconditionality assumption, e.g. if $\ell_1/X$
has a (UFDD) this can then be used to show that $\ell_1/X$ is the dual of
a subspace of $c_0.$

We would like to thank Gilles Godefroy and Dirk Werner for helpful
comments on the content of this note.

\section{Preliminary results}\label{prelim} \setcounter{equation}0

In this section we gather together some basic definitions and preliminary
results.

We start by recalling that a projection $P$ on a Banach space $X$ is an
{\it $L$-projection} if $\|x\|=\|Px\|+\|x-Px\|$ for any $x\in X.$  We
shall say that $P$ is a {\it $\theta-L$-projection} where $0<\theta\le 1$
if we have $\|x\|\ge \|Px\|+\theta\|x-Px\|.$  We shall say that $X$ is an
{\it L-summand} (respectively a
{\it $\theta-L-$summand}) if there is an
L-projection (respectively a $\theta-L-$projection of
$X^{**}$
onto $X;$ we shall say that $X$ is a {crude-$L$-summand} if it can be
equivalently renormed to be a $\theta-L-$summand for some $0<\theta\le
1.$ We also recall that
$X$ is called an
{\it
$M$-ideal} if
the canonical projection $\pi$ of $X^{***}$ onto $X^{*}$ is an
$L$-projection. Similarly
 that a Banach space $X$ is a {\it $\theta-M$-ideal} if $\pi$ is a
$\theta-L-$projection and a crude $M$-ideal if it has an equivalent
norm so that it is a $\theta-M-$ideal for some $0<\theta\le 1.$   For
background on the theory of M-ideals we refer to \cite{HWW}. The notion
of a crude M-ideal has also been considered in the literature,
originating with the work of Ando \cite{An} and more recently as a
special case of the so-called $M(r,s)-$inequalities, \cite {HO},
\cite{CN1} and \cite{CN2}.

Let us recall that $X$ has the {\it strong Schur property} \cite{R} if
there is a constant $c>0$ so that if $(x_n)$ is  any normalized sequence
 with
$\|x_m-x_n\|\ge \delta>0$ for any $m\neq n$ then there is a subsequence
$(x_{n})_{n\in \mathcal M}$ such that
$$ \|\sum_{k\in \mathcal M}\alpha_kx_k\|\ge
c\delta\sum_{k\in \mathcal M}|\alpha_k|$$  for any finitely non-zero
sequence $(\alpha_k)_{k\in \mathcal M}.$  This notion was first
introduced implicitly in Johnson-Odell \cite {JO}, and then explicitly by
Rosenthal \cite {R}
 and later studied by  Bourgain and Rosenthal \cite{BR}.

We will need some equivalent formulations of the strong Schur property:

\begin{Prop}\label{strong}  Let $X$ be a Banach space.  The following are
equivalent:
\newline
(i) $X$ has the strong Schur property.   \newline
(ii) There is a constant $c_1>0$ so that if  $(x_n)$ is any normalized
sequence with $\inf_{m>n}\|x_m-x_n\|=\delta$ then there exists $x^*\in
B_{X^*}$ with
$$\limsup_{n\to\infty}x^*(x_n)-\liminf_{n\to\infty}x_n^*(x)\ge
c_1\delta.$$\newline
(iii) For some fixed $\epsilon>0$ there exists a constant $c_2>0$ so that
if
$(x_n)$ is a normalized sequence
 with $\inf_{m>n}\|x_m-x_n\|\ge 1-\epsilon$ for any $m\neq n$ then
$x^*\in B_{X^*}$ with $\limsup_{n\to\infty} x^*(x_n)\ge c_2.$  \newline
(iv) There is a constant $c_3>0$ so that for any sequence $(x_n)$ in
$X$ there exists $x^*\in B_{X^*}$ with $\limsup_{n\to\infty} x^*(x_n)\ge
c_3\limsup_{n\to\infty}
\|x_n\|.$
\end{Prop}

\begin{proof} The equivalence of (i) and (ii) is essentially
contained in the usual proof of Rosenthal's $\ell_1-$theorem (cf.
\cite{D} pp. 209-211).  (ii) trivially implies (iii).

Let us now prove that (iii) implies (iv).  By the Uniform Boundedness
Principle we may suppose $(x_n)$ bounded and by passing to subsequences
and renormalizing we may suppose that $\|x_n\| =1$ for all $n.$  Let
$\delta=\inf_{m>n}\|x_m-x_n\|,$ and suppose $x^{**}$ is any
weak$^*$-cluster point of $(x_n).$ If
$\delta<1-\epsilon$ then
$\|x^{**}-x_n\|\le 1-\epsilon$ for all $n$ and so $\|x^{**}\|\ge
\epsilon.$  In this case there exists $x^*\in B_{X^*}$ with
$\limsup_{n\to\infty}x^*(x_n)\ge \frac12\epsilon.$  If $\delta\ge
1-\epsilon$, then we apply (ii).  Thus (iii) holds with
$c_3=\min(c_2,\frac{\epsilon}{2}).$

 Finally we show (iv) implies (ii). Indeed there exists $x^*\in B_{X^*}$
with $\limsup x^*(x_{2n}-x_{2n-1}) \ge c_3\delta$ and so
$$\limsup_{n\to\infty}x^*(x_n)-\liminf_{n\to\infty}x^*(x_n)\ge
c_3\delta.$$
\end{proof}

We will be interested in conditions which guarantee that a Banach space
$X$ into $c_0.$  We next state a criterion from \cite{KW} (the almost
isometric case) and \cite{GKL}.

\begin{Thm}\label{GKL}  Let $X$ be a separable Banach space.  Suppose
there is a constant $c>0$ so that if $x^*\in X^*$ and $(x^*_n)$ is any
weak$^*$-null sequence then
 $$\liminf_{n\to\infty}\|x^*+x_n^*\|\ge
\|x^*\|+c\liminf_{n\to\infty}\|x_n^*\|.$$
Then $X$ is isomorphic to a subspace of $c_0.$\end{Thm}

Note here that we can replace $\liminf$ by $\limsup$ or consider only the
case when both limits exist without changing the criterion.  A norm with
this property is called {\it Lipschitz-$UKK^*$}.  We now give a simple
application in the spirit of later results. We refer also to \cite{GKLi}
for connections between embeddability into $c_0$ and the strong Schur
property.

\begin{Thm}\label{Mideal}  Suppose $X$ is a separable Banach space such
that $X^*$ has the strong Schur property and suppose $X$ is a crude
$M$-ideal.  Then $X$ is isomorphic to a subspace of $c_0.$\end{Thm}

 \begin{proof} We may suppose $X$ is a $\theta-$M-ideal for some
$0<\theta\le 1.$  Suppose $x^*\in X^*$ and that $(x_n^*)$ is a
weak$^*$-null sequence.
Then there exists $x^{**}\in B_{X^{**}}$ so that $\limsup
x^{**}(x^*_n)\ge c_3\limsup\|x_n^*\|.$  Hence $(x_n^*)_{n=1}^{\infty}$
has a weak$^*$-cluster point $x^{***}\in X^{***}$ with $\|x^{***}\|\ge
c_3\limsup\|x_n^*\|.$  Clearly $x^{***}\in X^{\perp}$ and so
\begin{equation}\begin{align*}
 \limsup_{n\to\infty}\|x^*+x_n^*\| &\ge \|x^*+x^{***}\| \\
 &\ge \|x^*\|+\theta\|x^{***}\|\\
 &\ge \|x^*\| +c_3\theta\limsup_{n\to\infty}\|x_n^*\|.
\end{align*}\end{equation} We can now apply the result of \cite{GKL} to
deduce that $X$ embeds into $c_0.$\end{proof}

Another important concept we use concerns  unconditionality.
We shall say that a Banach space $X$ is of  {\it unconditional
type} if
whenever $x\in X$ and $(x_n)$ is a weakly null sequence in $X$ we have:
$$ \lim_{n\to\infty}(\|x+x_n\|-\|x-x_n\|)=0.$$  We shall say that $X$ is
of {\it shrinking unconditional type} if whenever $x^*\in X^*$ and
$(x_n^*)$ is weak$^*$-null in $X^*$ then
$$\lim_{n\to\infty}(\|x^*+x_n^*\|-\|x^*-x_n^*\|)=0.$$
These notions were introduced and studied (with different terminology) by
Neuwirth
\cite{N}.   We note first:

\begin{Lem} If $X$ is a separable Banach space which has shrinking
unconditional type then
$X$ has unconditional type.\end{Lem}

\begin{proof} Suppose $x\in X$ and $(x_n)$ is weakly null and that
$\|x+x_n\|>\|x-x_n\|+\epsilon$ for all $n,$ where $\epsilon>0.$  Choose
$y_n^*\in B_{X^*}$ so that $y_n^*(x+x_n)=\|x+x_n\|.$  By passing to a
subsequence  we can suppose $y_n^*$ converges to some $x^*\in X^*.$
Then
$\lim_{n\to\infty}\|2x^*-y_n^*\|=\lim_{n\to\infty}\|y_n^*\|=1.$  Now
$\lim_{n\to\infty} (\|x+x_n\|-y_n^*(x_n))=x^*(x)$ and so $$
\lim_{n\to\infty}(\langle
x-x_n,2x^*-y_n^*\rangle-\|x+x_n\|) = 0.$$  This implies that
$
\liminf(\|x-x_n\|-\|x+x_n\|)\ge 0$ and gives the lemma.\end{proof}

Let us recall that a separable Banach space $X$ has (UMAP) if there is a
sequence of finite-rank operators $(T_n)$ such that
$\lim_{n\to\infty}T_nx=x$ for $x\in X$ and
$\lim_{n\to\infty}\|I-2T_n\|=1$ (see \cite{CK}, \cite{GK}); we say $X$
has {\it shrinking (UMAP)} if, in addition
$\lim_{n\to\infty}T_n^*x^*=x^*$ for $x^*\in X^*.$
 It is shown in
\cite{GK}
that $X$ has (UMAP) if and only if for every $\epsilon>0$ $X$ is
isometric to a one-complemented subspace of a space $V_{\epsilon}$ with a
$(1+\epsilon)$-(UFDD).

\begin{Lem} Let $X$ be a Banach space with (UMAP); then $X$ is
of unconditional type. If $X$ has shrinking (UMAP) $X$ is of
shrinking unconditional type. \end{Lem}

\begin{proof} Suppose $x\in X$ and $(x_n)$ is weakly null. It is enough
to show that $\lim_{n\to\infty}\|x+x_n\|\le \lim_{n\to\infty}\|x-x_n\|$
under the assumption that both limits exist.
   \begin{equation}\begin{align*}
\lim_{n\to\infty}\|x+x_n\|&=\lim_{k\to\infty}\lim_{n\to\infty}\|(2T_k-1)x+x_n\|
\\
&=\lim_{k\to\infty}\lim_{n\to\infty}\|(2T_k-I)x+(I-2T_k)x_n\|\\
&\le \lim_{n\to\infty}\|x-x_n\|.\end{align*}\end{equation}
The shrinking case is similar. \end{proof}

\begin{Lem} Let $X$ be a separable Banach space of shrinking
unconditional type.
Then any subspace or quotient of $X$ has shrinking unconditional type.
\end{Lem}

\begin{proof} If $Y$ is a subspace of $X$ then $(X/Y)^*$ can be
identified with $Y^{\perp}$ and trivially $X/Y$ has  shrinking
unconditional type.  Now $Y^*$ can be identified with $X^*/Y^{\perp}.$
Let $Q:X^*\to Y^*$ be the canonical quotient map.  Suppose $y^*\in Y^*$
and that $(y_n^*)$ is weak$^*$-null in $Y^*$.  Suppose that
$\|y^*+y_n^*\|<\|y^*-y_n^*\|-\epsilon$ where $\epsilon>0.$  We may pick
(by the Hahn-Banach theorem)
$x_n^*\in X^*$ so that $\|x_n^*\|=\|y^*+y_n^*\|$ and $Qx_n^*=y^*+y_n^*.$
Passing to a subsequence we can suppose that $x_n^*$ converges weak$^*$
to $x^*.$  Now
$$ \lim_{n\to\infty}(\|2x^*-x_n^*\|-\|x_n^*\|)=0$$
and $Q(2x^*-x_n^*)=2Qx^*-y^*-y_n^*= y^*-y_n^*$ by the weak$^*$-continuity
of $Q.$  Hence
$$ \limsup_{n\to\infty}(\|y^*-y_n^*\|-\|y^*+y_n^*\|) \le 0,$$
which yields a contradiction, so that $Y$ is of shrinking unconditional
type. \end{proof}

\begin{Lem}\label{unctype}  Let $X$ be a subspace of a space with (UMAP);
if
$X$ does not contain $\ell_1$ then $X$ has shrinking unconditional
type.\end{Lem}

\begin{proof} Suppose $\epsilon>0$ and that $x^*\in X^*$ and that
$(x_n^*)$ is any weakly null sequence; assume that
$\sup_n\|x^*+x_n^*\|\le 1$ and
$\|x^*-x_n^*\|>\|x^*+x_n^*\|+\epsilon$ for all $n\in\mathbb N$
where $\epsilon>0.$  By results of
\cite{GK} we may suppose
$X$ is isometric to a subspace of a space $V=V_{\epsilon}$ with a
$(1+\epsilon/2)-$(UFDD).  Indeed suppose $(Q_n)$ are finite rank
projections defining a $(1+\epsilon/2)-$(UFDD).  Let
$T_n=\sum_{k=1}^nQ_i.$ Let
$j:X\to V$ be the isometric embedding.

If $v^*\in V^*$ we have $j^*v^*=\sum_{k=1}^{\infty}j^*Q_k^*v^*$
unconditionally in the
weak$^*$-topology.  Since $X^*$ does not contain $c_0$ (or equivalently
$\ell_{\infty}$) this series converges in norm so that
$\lim_{k\to\infty}\|j^*v^*-j^*T_k^*v^*\|=0.$

Now by the Hahn-Banach theorem we can find $v_n^*\in V^*$ so that
$\|v_n^*\|=\|x^*+x_n^*\|$ and $j^*v_n^*=x^*+x_n^*.$   By passing to a
subsequence we can suppose that $v_n^*$ converges weak$^*$ to some $v^*.$
Clearly $j^*v^*=x^*.$ Then
\begin{equation}\begin{align*}
 \limsup_{n\to\infty}\|x^*-x_n^*\| &=
\limsup_{n\to\infty}\|2j^*v^*-j^*v_n^*\|\\
&=
\limsup_{k\to\infty}\limsup_{n\to\infty}\|2j^*T_k^*v^*-j^*v_n^*-2j^*(T_k^*v^*-T_kv_n^*)\|
\\
&\le\limsup_{k\to\infty}\limsup_{n\to\infty}\|(2T_k-I)v_n^*\|\\
&\le
(1+\frac{\epsilon}{2})\limsup_{n\to\infty}\|x^*+x_n^*\|.\end{align*}\end{equation}
This contradiction establishes the lemma.\end{proof}

\section{Subspaces of $c_0$ and trees}\label{trees}
\setcounter{equation}0

Consider the set ${\cal F}{\Bbb N}$ of all finite subsets of $\Bbb N$
with the following partial order.  If $a=\{n_1,n_2,\ldots,n_k\}$ where
$n_1<n_2<\ldots<n_k$ and $b=\{m_1,m_2\ldots,m_l\}$ where
$m_1<m_2<\cdots<m_l,$ then $a\le b$ if and only if $k\le l$ and
$m_i=n_i$ where $1\le i\le k$ (i.e. $a$ is an initial segment of $b.$)
We say that $b$ is a {\it successor} of $a$ if $|b|=|a|+1$ and $a\le b;$
the collection of successors of $a$ is denoted by $a+$.
If $a\neq\emptyset$ then $a-$ denotes the unique predecessor of $a$ i.e.
$a$ is a successor of $a-.$ Let
$S$ be a subset of
$\cal F\Bbb N.$ We will say that
$S$ is a
{\it full tree}
if we have
 \begin{enumerate} \item $\emptyset \in S.$
  \item  Each
$a\in S$ has infinitely many
successors in $S.$
\item If $a\in S$ and $\emptyset\neq a\in S$ then
$a-\in S$.
\end{enumerate}
  It is easy to see that any full tree
 is isomorphic as an ordered set to $\cal F\Bbb N.$
 If $S$ is any full tree we will say that a sequence
$\beta=\{a_n\}_{n=0}^{\infty}$ is a {\it branch} of $S$ if $a_n\in S$
for all
$n,$ $a_0=\emptyset$ and $a_{n+1}$ is a successor of $a_n$ for all
$n\ge
0.$

Now let $V$ be a vector space. We define a {\it tree-assignment} to be a
map $a\to x_a$ defined on a full tree $S$. We define a
{\it tree-map} to be a tree-assignment
$a\to x_a$  with the properties that
$x_{\emptyset}=0$ and for every branch $\beta$ the set
$\{a:\ x_{a}\neq 0:\ a\in\beta\}$ is finite.  Given any tree-map we
define a {\it height function} $h$ which assigns to each $a$ a countable
ordinal; to do this we define $h(a)=0$ if $x_{b}=0$ for $b\ge a$ and
then inductively $h(a)$ is defined by $h(a)\le \eta$ if and only if
$h(b)<\eta$ for every $b>a.$  The {\it height} of the tree-map is
defined to be
$h(\emptyset).$  Note that the tree-map $a\to x_{a}$ has finite height
$m\le n$ if and only if $x_{a}=0$ whenever $|a|>n.$

The following easy lemma, proved in \cite{GKL2} is a restatement of the
fact that certain
types of games (which are not used in this paper) are determined.

\begin{Lem}\label{determine}Suppose $(x_a)_{a\in S}$ is a tree-map and
that $A$ is any subset of $V.$  Then either there is a full tree
$T\subset S$ so that $\sum_{a\in\beta}x_a\in A$ for every branch
$\beta\subset T$ or there is a full tree $T\subset S$ so that
$\sum_{a\in\beta}x_a\notin A$ for every branch $\beta\subset
T.$\end{Lem}

Now suppose $V=X$ is a Banach space.
If $\tau$
is a topology on $X$ (e.g. the weak topology or for dual spaces the
weak$^*$-topology) we say that a tree-map $(x_a)_{a\in S}$ is {\it
$\tau-$null} if for every $a\in S$ the set $\{x_b\}_{b\in a+}$ is a
$\tau$-null sequence.

We now introduce a definition which will characterize subspaces of $c_0.$
We say that a Banach space $X$ has the {\it bounded tree property} with
constant $\sigma>0$ if every weakly null tree-map $(x_a)_{a\in S}$  has
a full subtree $T$ so that $\|\sum_{a\in\beta}x_a\|\le 1$ for every
branch $\beta.$

\begin{Thm}\label{boundedtree} Let $X$ be a separable Banach space
containing no copy of
$\ell_1$ with the bounded tree property.  Then $X$ is isomorphic to a
subspace of $c_0.$\end{Thm}

\begin{proof} (Compare Theorem 4.1 of \cite{GKL2}.)  Define for $x\in X,$
$f(x)$ to be the infimum of all $\lambda>0$ so that for every weakly
null tree-map $(x_a)_{a\in S}$ with $\|x_a\|\le \sigma$ there is a full
tree $T\subset S$ with $\|x+\sum_{a\in\beta}x_a\|\le \lambda$ for every
branch $\beta.$  Note that $\|x\|\le f(x)$, $f(0)\le 1$, $f(-x)=f(x)$ and
that
$|f(x)-f(y)|\le \|x-y\|.$ In particular we have $\|x\|\le f(x)\le
\|x\|+1.$ We now argue exactly as in
\cite{GKL2} that
$f$ is convex.  For the convenience of the reader we repeat the
argument.
Let $u=t x+(1-t)y,$ where $0<t<1.$
Suppose
$\lambda>f(x)$ and $\mu>f(y).$  Let $(x_a)_{a\in S}$ be any weakly
null tree-map of height $k$ with $\|x_a\|\le \sigma$ for all $a\in S.$
Then we can find a full subtree $T_1\subset S$ so that for every branch
$\beta$ we have
$$ \|x+\sum_{a\in\beta}x_a\|\le \lambda$$
and then a full subtree $T_2\subset T_1$ so that for every branch
$\beta\subset T_2$
$$ \|y+\sum_{a\in\beta}x_a\| \le \mu.$$
Obviously for every branch $\beta\subset T_2$
$$ \|u+\sum_{a\in\beta}x_a\| \le t\lambda+(1-t)\mu$$
so that $f_k(u)\le t\lambda+(1-t)\mu.$

Next we note that if $\|x_n\|\le \sigma$ and $\lim_{n\to\infty}x_n=0$
weakly then $\limsup f(x+x_n)\le f(x).$
Assume that
$\lambda <\limsup_{n\to\infty}f(x+x_n).$  By passing to a
subsequence we can suppose $\lambda<f(x+x_n)$ for every $n.$ Then for
each
$n$ there is a weakly null tree-map  $(y^{(n)}_a)_{a\in S_n}$ of height
$k$ so that $\|y^{(n)}_a\|\le \sigma$ for all $a\in S_n$ and $$
\|x+x_n+\sum_{a\in\beta}y^{(n)}_a\| >\lambda$$ for every branch
$\beta\subset S_n.$
Now let $T$ be the tree consisting of all sets $\{m_1,\ldots,m_l\}$
where
$m_1<m_2<\cdots<m_l$ such that if $l>1$ then $\{m_2,\ldots,m_l\}\in
S_{m_1}.$  We define a weakly null tree-map  by
$$ z_{m_1,\ldots,m_l}=\begin{cases} x_{m_1}  \text{ if } l=1\\
y^{(m_1)}_{m_2,\ldots,m_l} \text{ if } l>1.\end{cases}$$
Then for every branch $\beta\subset T$  we have
$$ \|x+\sum_{a\in\beta}z_a\| >\lambda$$
so that $f
(x)\ge \lambda$.
This implies our claim.

Let $|\cdot|$ be the Minkowski functional of the set $\{x:f(x)\le 2\}$.
Then $|\cdot|$ is a norm on $X$ satisfying $\frac12\|x\|\le |x|\le
\|x\|.$  Suppose $|x|=1$ and $(x_n)$ is weakly null sequence with
$|x_n|\le \frac12\sigma.$  Then $\|x_n\|\le \sigma$ and so
$ \limsup f(x+x_n)\le 2.$    Hence
$$ \limsup|x+x_n|\le 1.$$

Now by Proposition 2.7 of \cite{GKL2} we have that $X^*$ is separable.
We  can  then apply Proposition 2.6 of \cite{GKL2} to deduce that if
$|x^*|=1$ and $(x_n^*)$ is a weak$^*$-null sequence in $X^*$ with
$\|x_n^*\|=\tau$ then
$$ \liminf_{n\to\infty}|x^*+x_n^*|\ge 1+\frac{\sigma}{12}\tau.$$
Thus $X^*$ has a Lipschitz-UKK$^*$ norm and by the results of \cite{GKL}
(see also \cite{KW}) this implies that $X$ embeds into $c_0.$
\end{proof}

We next introduce a dual notion.  We say that $X^*$ has the {\it
weak$^*$ summable tree property } with constant $c>0$ if
 for every weak$^*$-null tree-map $(x_a^*)_{a\in S}$ on $X^*$
satisfying the boundedness property: \begin{equation}\label{bp}
\sup_{a\in S}\|\sum_{b\le a}x_b\| <\infty, \end{equation} and  for every
$\epsilon>0$  there is a full subtree $T$ so that
$$\|\sum_{a\in\beta}x_a^*\| >c\sum_{a\in\beta}\|x_a^*\|-\epsilon$$ for
every branch $\beta.$ Notice that if $X$ is a subspace of $c_0$ then
\cite{KW}
$$ \liminf_{n\to\infty}\|x^*+x_n^*\| \ge \|x^*\|
+\liminf_{n\to\infty}\|x_n^*\|$$ and this implies directly that $X^*$
has the weak$^*$ summable tree property with constant one.

\begin{Thm}\label{summabletree} Suppose $X$ is a separable Banach space
such that $X^*$ has the weak$^*$-summable tree property.  Then $X$ is
isomorphic to a subspace of $c_0.$\end{Thm}

\begin{proof} We show that $X$ contains no subspace isomorphic to
$\ell_1$ and that $X$ has the bounded tree property.  To show $X$
contains no copy of $\ell_1$ it suffices to show that $\ell_2$ does not
embed in $X^*$ by \cite{P}.  Suppose then $u_n^*$ is a weak$^*$-null
sequence in $X^*$ so that $\|u_n^*\|=1$ and $\|\sum_{i\in A}u_i\|\le
C|A|^{1/2}$ for any finite subset of $A$ of $\mathbb N$ where $C$ is an
absolute constant.  Then define, for any $N,$ the tree-map on $\cal
F\mathbb N$ by
        $$ x_a^*=u_{m_n}^* \qquad \text{if }a=\{m_1,\ldots,m_n\}\text {
and }1\le n\le N$$
and $x_a^*$ otherwise.  It is clear that any full subtree $T$ has a
branch $\beta$ with $\|\sum_{a\in\beta}u_a^*\|\le CN^{1/2}$ while
$\sum_{a\in\beta}\|u_a^*\|=N$, and therefore if $N$ is large enough we
obtain a contradiction to the weak$^*$ summable tree property.

 Next we need a duality argument.  We assume that $X^*$ has the weak$^*$
summable tree property with constant $c>0.$  We show that $X$ has the
bounded tree property with constant $\sigma$ for any $0<\sigma<c/2.$
Indeed if not there is by Lemma \ref{determine} a weakly null tree-map
$(x_a)_{a\in S}$ with the properties that $\|x_a\|\le \sigma$ for all
$a$ and $\|\sum_{a\in \beta}x_a\|>1$ for every branch $\beta.$   For
each branch $\beta$ pick $u_{\beta}^*\in X^*$ with
$\|u_{\beta}^*\|=1$ and $\langle
\sum_{a\in\beta}x_a,u_{\beta}^*\rangle >1.$ Let
$h$
be the height function of the given tree-map.
For each $a\in S$ we define $y_a^*$ by transfinite induction on $h(a).$
If $h(a)=0$ let $y_a^*=u_{\beta}^*$ where $\beta$ is any branch to which
$a$ belongs.  Then if $(y_a^*)$ has been defined for $h(a)<\eta$ and if
$h(b)=\eta$ we define $(y_b^*)$ to be any weak$^*$-cluster point of
$(y_a^*)_{a\in b+}.$  (Note that according to our definition $(y_a^*)$
is a tree-assignment but not necessarily a tree-map because it is not
supported
on a well-founded tree and we may have $y^*_{\emptyset}\neq 0.$)

Let us now make a tree-map by defining $x^*_{\emptyset}=0$ and then if
$h(a-)\ge 1$ we define $x^*_a=y^*_a-y^*_{a-}.$  If $h(a-)=0$ we define
$x_a^*=0.$  This is clearly a tree-map which also satisfies (\ref{bp})
and we have that for each
$a\in
S,$ zero is a weak$^*$-cluster point of $(x^*_b)_{b\in a+}.$  It is then
easy to see that we can pass to a full subtree $T$ so that
$(x^*_a)_{a\in T}$ is weak$^*$-null.  Let $x^*=y^*_{\emptyset}.$

Now pick  $\epsilon>0$  so that $3\epsilon+2c^{-1}\sigma<1.$ We can use
the definition of the weak$^*$-summable
tree property and also Lemma 3.3 of \cite{GKL2} to pass to a further
full subtree (still labeled $T$) so that we have $|\langle
x_a,x^*\rangle|< \epsilon/2^{|a|}$ when $|a|>0$ and for any branch
$\beta\subset T$
\begin{equation}\begin{align*}
|\langle \sum_{a\in\beta}x_a,\sum_{a\in\beta}x_a^*\rangle
-\sum_{a\in\beta}\langle x_a,x_a^*\rangle|&\le \epsilon\qquad\\
c(\sum_{a\in\beta}\|x_a^*\|-\epsilon)&\le
 \|\sum_{a\in\beta}x_a^*\|.\qquad
\end{align*}\end{equation}

For any branch $\beta$ let $b$ be the first point where $h(b)=0$.  Then
$$\|x^*+\sum_{a\in\beta}x_a^*\|=\|y^*_b\|=1.$$  It follows that
$$ \|\sum_{a\in\beta}x_a^*\| \le 2.$$
Now  we have
\begin{equation}\begin{align*}
1 &< \langle \sum_{a\in\beta}x_a,y_b^*\rangle\\
&\le \sum_{a\in\beta}|\langle x_a,x^*\rangle| +|\langle \sum_{a\in\beta}
x_a, \sum_{a\in\beta}x_a^*\rangle|\\
&\le 2\epsilon +\sum_{a\in\beta}|\langle x_a,x_a^*\rangle|\\
&\le  2\epsilon +\sigma \sum_{a\in\beta}\|x_a^*\|\\
&\le 3\epsilon + c^{-1}\sigma \|\sum_{a\in\beta}x_a^*\|\\
&\le 3\epsilon + 2c^{-1}\sigma.\end{align*}\end{equation}
This gives a contradiction and so we deduce that $X$ has the bounded
tree property and we can apply Theorem \ref{boundedtree} to obtain the
result.\end{proof}

\section{The very strong Schur
property}\label{vstrongschur}\setcounter{equation}0

We shall say that a tree-assignment $(x_a)_{a\in S}$ in $X$ is {\it
$\delta$-separated} if $\|x_{b}-x_{b'}\|\ge \delta$ whenever $b,b'\in S$
are such that $b,b'\in a+$ for some $a\in S.$ Let us say that a Banach
space
$X$ has the
{\it very strong Schur
property} if there is a constant $c>0$ so that whenever $(x_a)_{a\in S}$
is a $\delta$-separated bounded tree-assignment then there is branch
$\beta$ and
$x^*\in B_{X^*}$ with $|x^*(x_a)|\ge c\delta$ whenever $\emptyset\neq
a\in
\beta.$

We first justify this terminology:

\begin{Prop}\label{verystrong} Suppose $X$ is a Banach space with the
very
strong Schur property.  Then $X$ has the strong Schur property.\end{Prop}

\begin{proof} We verify condition (iii) of Proposition \ref{strong} with
$\epsilon=\frac12.$ Let $(x_n)$ be a normalized sequence with
$\inf_{m>n}\|x_m-x_n\|\ge \frac12.$  Form a tree-assignment
$(y_a)_{a\in{\mathcal F}{\mathbb N}}$ by putting $y_{\emptyset}=0$ and
then $y_{a}=x_{n}$ if $n=\max a.$  Then $(y_a)$ is a bounded
$\frac12-$separated tree-assignment and so there is a branch $\beta$  and
$x^*\in B_{X^*}$ with $|x^*(y_a)|\ge \frac12c.$  This leads to a
subsequence $(x_{n_k})_{k=1}^{\infty}$ where $|x^*(x_{n_k})|\ge \frac12c$
and Proposition \ref{strong} (iii) holds with either $x^*$ or $-x^*.$
\end{proof}

There is an important situation when the converse is true.

\begin{Thm}\label{verystrong2}  Suppose that $Y$ is a crude
$L$-ideal.  If $X$ is a closed subspace of $Y$ with the strong Schur
property then $X$ has the very strong Schur property.\end{Thm}

\begin{proof} We may suppose that $Y$ is a $\theta$-L-ideal where
$0<\theta\le 1.$  Suppose $P$ is the associated $L$-projection.  We also
use Proposition \ref{strong} (iv) to deduce there is a constant $c>0$ so
that if
$(x_n)_{n\in\mathbb N}$
is any bounded sequence in $X$ with $\inf_{m\neq n}\|x_m-x_n\|\ge
\delta>0$ then
$(x_n)$ has a subsequence $(w_n)$ such that we have an estimate
\begin{equation}\label{ell1} \|\sum_{k=1}^{\infty}\alpha_kw_k\| \ge
c\sum_{k=1}|\alpha_k|\end{equation} for all finitely nonzero sequences
$(\alpha_k).$

Now suppose $(x_a)_{a\in S}$ is a $\delta-$separated tree-assignment.
Let $\sigma=\frac14c\theta$.  We shall show by an inductive
construction that there is a branch
$\beta$ and for each $a\in\beta,$ $x_a^*\in X^*$ with $\|x_a^*\|<1$ so
that
$|x_a^*(x_b)|\ge \sigma\delta$ if $\emptyset \neq b\le a.$  This will
complete the proof since then we can take $x^*$ as any weak$^*$ cluster
point of $\{x_a^*:a\in\beta\}.$

We start the branch with $\emptyset.$  Now suppose $a\in\beta$; we must
choose a successor $b\in a+$ and a corresponding $x_b^*.$ First let
$y^*_a$ be any norm-preserving extension of $x_a^*$ to $Y.$  Next we pick
a subsequence $w_n=x_{b_n}$ of $\{x_b: \  b\in a+\}$ satisfying
(\ref{ell1}).  Let $x^{**}$ be any weak$^*$-cluster point of
$(w_n)_{n=1}^{\infty}.$

Suppose $y\in Y.$  Let $x^{**}-y$ belongs to the weak$^*$-closed convex
hull
$W_k$ of
$\{w_n-y\}_{n=k}^{\infty}$ and hence $0$ is in the norm-closure of the
set
$W_k-\|x^{**}-y\|B_Y.$  We deduce that for any $\epsilon>0$ we can find
convex combinations $\sum_{j=1}^k\alpha_j(w_n-y)$ and
$\sum_{j=k+1}^l\alpha_j(w_n-y)$ of norm at most
$\|x^{**}-y\|+\epsilon/2.$  Hence
$$ 2c\delta \le \|\sum_{j=1}^k\alpha_jw_j-\sum_{j=k+1}^l\alpha_jw_j\| \le
2\|x^{**}-y\|+\epsilon.$$
Thus $d(x^{**},Y) \ge c\delta.$

In particular, $\|x^{**}-Px^{**}\| \ge c\delta.$  Let $E$ be the linear
span of $\{x_d:d\le a\}\cup\{Px^{**},x^{**}-Px^{**}\}.$  We define a
linear functional $\varphi$ on $E$ by
$\varphi(e)=y_a^*(e)$ if $e\in E\cap Y$ and
$\varphi(x^{**}-Px^{**})=2\sigma$ if $y_a^*(Px^{**})\ge 0$
and
$\varphi(x^{**}-Px^{**})=-2\sigma$ if $y_a^*(Px^{**})<0.$  For any $e\in
E$ we have  $e=e_0+\lambda (x^{**}-Px^{**})$ where $e_0\in E\cap Y$ and
$\lambda\in\mathbb R$.  Then
$$ |\varphi (e)|\le |y_a^*(e_0)|+2|\lambda||\sigma|  \le \|x_a^*\|\|e_0\|
+\frac12 \theta\|e-Pe\|.$$
Hence $\|\varphi\|<1.$  It follows that $\varphi$ has a
weak$^*$-continuous extension $y^*\in Y^*$ with $\|y^*\|<1.$
Now $|\langle y^*,x^{**}\rangle| \ge 2\sigma$ and hence we can pick $n$
so that $|y^*(w_n)|=|y^*(x_{b_n})|\ge \sigma.$  We thus select $b=b_n$
and set $x_b^*=y^*|_X.$  This inductive process establishes our
result.\end{proof}

Let us now remark that a closed subspace of $L_1$ has the very strong
Schur property if and only if it has the strong Schur property since
$L_1$ is an $L$-ideal in its bidual.  It follows therefore that the
examples constructed by Bourgain and Rosenthal \cite{BR} show that
for subspaces of $L_1,$ the very strong Schur property does not imply
embeddability into $\ell_1$ or even the Radon-Nikodym Property.  However,
Johnson and Odell
\cite{JO} showed
that a subspace of $L_1$ with a (UFDD) and the strong Schur property is
isomorphic to a subspace of $\ell_1;$ see also \cite{O}. Thus the
presence of
some unconditionality is crucial here.  This motivates our next theorem:

\begin{Thm}\label{c0dual} Let $X$ be a separable Banach space with the
property that $X^*$ has the very strong Schur property.  Assume that $X$
is linearly isomorphic to a subspace of a Banach space with (UFDD).  Then
$X$ is linearly isomorphic to a subspace of $c_0.$\end{Thm}

\noindent{\it Remark.} The assumption that $X$ embeds into a
space with
(UFDD) is equivalent to the assumption that $X$ embeds in a space with
unconditional basis \cite{LT} p. 51.  As will be seen in the proof
the theorem holds if $X$ is assumed to have shrinking
unconditional type.

\begin{proof} Note first that $X$ cannot contain $\ell_1$ by results of
\cite{P} since $X^*$ has the Schur property.  Therefore we can apply
Lemma \ref{unctype} to deduce that $X$ can be given an equivalent norm
so that it has shrinking unconditional type.  We complete the proof by
showing that $X$ has the weak$^*$-summable tree  property and applying
Theorem \ref{summabletree}.

Assume that $X^*$ has the very strong Schur property with constant $c.$
We will show that $X$ has the weak$^*$-summable tree property with
constant $c/2.$  Suppose $(x_a^*)_{a\in S}$ is a weak$^*$-null tree such
that
$$ \sup_{a\in S}\|\sum_{b\le a}x_b^*\|=M<\infty.$$
Assume that $(x_a^*)$ fails to have a full subtree such that
$$ \|\sum_{a\in\beta}x_a^*\|
>\frac{c}{2}\sum_{a\in\beta}\|x_a^*\|-\epsilon.$$
Then by considering the tree-map $(x_a^*,\|x_a^*\|)$ in $X^*\times
\mathbb R$ and using Lemma \ref{determine} we can find a full subtree
$(x_a^*)_{a\in S_1}$ so that for every branch we have:
$$ \|\sum_{a\in\beta}x_a^*\| \le \frac{c}{2} \sum_{a\in
\beta}\|x_a^*\|-\epsilon.$$

Next we can pass to a full subtree $S_2$ so that if for each $a\in S_2$
either $\inf_{b\in a+}\|x^*_a\|>0$ or $\sup_{b\in a+}\|x^*_b\|\le
2^{-|a|-3}\epsilon.$
We then define a tree assignment $(u_b^*)$ as follows.  Put
$u_{\emptyset}^*=0.$  If $a\in S_2$ is such that $\inf_{b\in
a+}\|x_b^*\|=0$ then let $\{u_b^*:b\in a+\}$ be assigned to be any fixed
weak$^*$-null normalized sequence.  If $a\in S_2$ and $\inf_{b\in
a+}\|x^*_b\|>0$   we let
$u^*_b=x^*_b/\|x^*_b\|$ if
$b\in a+.$
Then $(u_a^*)_{a\in S_2}$ is weak$^*$-null and, using the
weak$^*$-lower-semicontinuity of the norm we may pass to a full subtree
$S_3$ so that for any $a\in S_3$ we have
$$ \inf_{b,b'\in a+}\|u^*_b-u^*_{b'}\|\ge \frac12.$$

Next we use the fact that $X$ has shrinking unconditional type.  For each
$a\in S_3$ there is a closed absolutely convex weak$^*$-neighborhood of
the origin $W_a$ so that if $w^*\in W_a$ and $\|w^*\|\le 2M$ then
$$ \left|\sum_{b\le a}\eta_bx^*_b+w^*\|-\|\sum_{b\le a}\eta_bx_b^*-w^*\|
\right|
\le 2^{-|a|-2}\epsilon$$
for every choice of signs $\eta_b=\pm1$ for $b\le a.$
Let $T=\{a\in S_3:\ x^*_a \in 2^{|b|-|a|}W_b,\ \text{if }b<a\}.$
Then $T$ is a full subtree of $S_3.$

Let $\beta$ be any branch in $T.$  We write
$\beta=\{a_0,a_1,a_2,\ldots,\}$ where $a_0<a_1<a_2\ldots.$  Let
$$
\sigma_n:=\max_{\eta_k=\pm1}\|\sum_{k=0}^n\eta_kx_{a_k}^*
+\sum_{k=n+1}^{\infty} x_{a_k}^*\|.$$  Then
$\sigma_0=\|\sum_{k=1}^{\infty}x^*_{a_k}\|.$  Notice that if
$\eta_{n}=-1$ where $n\ge 1$ then since
$\sum_{k=n+1}^{\infty}x_{a_k}^*\in W_{a_n}$ and
$\|\sum_{k=n+1}^{\infty}x^*_{a_k}\|\le 2M,$
$$ \|\sum_{k=0}^n\eta_kx_{a_k}^*+\sum_{k=n+1}^{\infty}x_{a_k}^*\| \le
2^{-n-2}\epsilon +
\|\sum_{k=0}^{n-1}\eta_kx_{a_k}^*-\sum_{k=n}x_{a_k}^*\|.$$
Thus
$$ \sigma_n \le \sigma_{n-1} + 2^{-n-2}\epsilon.$$
It follows that
$$ \sigma_n \le \sum_{k=1}^n 2^{-k-2}\epsilon +
\|\sum_{k=0}^{\infty}x_{a_k}^*\|.$$
Thus we conclude that for any branch and any choice of signs $\eta_a$ we
have
$$ \|\sum_{a\in \beta}\eta_ax_a^*\| \le \|\sum_{a\in\beta}x_a^*\| +
\frac14\epsilon.$$

Next we can use the very strong Schur
property and the fact that $(u^*_a)_{a\in T}$ is $\frac12-$separated to
find a branch $\beta$ and $u^{**}\in B_{X^{**}}$ with $|u^{**}(u_a^*)|\ge
\frac{c}{2}$ for $a\in\beta.$  By the construction of $(u^*_a)$ we have
$$ \|x_a^*- \|x_a^*\|u_a^*\| \le 2^{-|a|-2}\epsilon$$
so that
$$ |u^{**}(x_a^*)| \ge \frac{c}{2}\|x_a^*\| - 2^{-|a|-2}\epsilon.$$
Choose $\eta_a=\pm 1$ so that $\eta_au^{**}(x_a^*)\ge 0.$  Then
$$ u^{**}(\sum_{a\in\beta}\eta_ax_a^*) \ge \frac{c}{2}
\sum_{a\in\beta}\|x_a^*\| - \frac12\epsilon.$$
Hence
$$ \|\sum_{a\in\beta}x_a^*\| \ge \frac{c}{2}\sum_{a\in
\beta}\|x_a^*\|-\frac34 \epsilon.$$
This is a contradiction and shows that $X$ has the weak$^*$ summable tree
property.  The proof is complete.\end{proof}

\begin{Cor} Suppose $X$ and $Y$ are Banach spaces such that $X^*$
and $Y^*$ are isomorphic.  Suppose $X$ is isomorphic to a subspace of
$c_0$.  Then $Y$ is isomorphic to a subspace of $c_0$ if and only if $Y$
embeds in a space with (UFDD).\end{Cor}

\noindent{\it Remark.} We do not know if one can conclude that $X$ and
$Y$ are isomorphic if both embed into $c_0.$

\section{The extension property}

Let us recall that if $X$ is a Banach space $X$ and $E$ is a closed
subspace of $X$ then the pair $(E,X)$ is said to have the {\it
$\lambda-$extension
property, $\lambda-$(EP)} if, for any compact Hausdorff space $K,$ every
bounded operator
$T:E\to C(K)$ has a bounded extension $\tilde T: X\to C(K)$ with
$\|\tilde T\|\le \lambda \|T\|$
(Johnson-Zippin \cite{JZ}).  We say $(E,X)$ has (EP) if it has
$\lambda-$(EP) for some $\lambda\ge 1.$ Johnson and Zippin
\cite{JZ} showed that if
$X$
is weak$^*$-closed subspace of $\ell_1=c_0^*$ then $(X,\ell_1)$ has (EP),
although curiously it is unknown whether it has $(1+\epsilon)-(EP)$ for
any $\epsilon>0.$  See \cite{S} and \cite{Z} for recent progress on
extension properties.

As observed in Corollary 1.1 of \cite{JZ}, using the results of
\cite{LR}, the extension property of
$(X,\ell_1)$ depends only on the quotient space $\ell_1/X$; hence it
follows that if $\ell_1/X$ is isomorphic to $Y^*$ where $Y$ is a
closed subspace of $c_0$ then $(X,\ell_1)$ has the extension property
(because there is an automorphism $\tau$ of $\ell_1$ so that $\tau(X)$ is
weak$^*$-closed).  The aim of this section is to show how the results of
the paper can give a partial converse to this theorem.

\begin{Thm}\label{JZSchur}  Suppose $X$ is a closed subspace of $\ell_1$
so that $(X,\ell_1)$ has (EP).  Then $\ell_1/X$ has the very strong Schur
property.\end{Thm}

\noindent{\it Remark.}  The result that $\ell_1/X$ has the Schur property
was obtained earlier by the author and A. Pe\l czy\'nski by somewhat
similar arguments.  This answered a question of Zippin concerning the
case $\ell_1/X\approx L_1.$

\begin{proof} We suppose that $(X,\ell_1)$ has $\lambda-$(EP). Let
$Y=\ell_1/X$ and denote by $Q_Y$ the quotient map of $\ell_1$ onto $Y$.

We start by supposing that $(y_a)_{a\in S}$ is a bounded
$\delta-$separated
tree assignment in $Y=\ell_1/X.$  Let $E_n$ be an increasing sequence of
finite-dimensional subspaces of $Y$ whose union is dense.
We start by observing that for each $a\in S$ and each $n\in\mathbb N$
there is an infinite number of $b\in a+$ so that $d(y_b,E_n)>
\delta/4.$ Indeed, if not there are infinitely many $b\in a+$ so that
$d(y_b,E_n)\le \delta/4$ and for each such $b$ we can find $e_b\in E_n$
with $\|y_b-e_b\|\le \delta/4.$  The set of such $e_b$ is bounded and so
by compactness arguments we obtain $b\neq b'$ with $\|y_b-y_{b'}\|\le
3\delta/4.$

Now we may pass to a full subtree $(y_a)_{a\in T}$ so that there exists a
map $\psi:T\setminus\{\emptyset\}\to \mathbb N$ with the properties
that $d(y_a,E_{\psi(a)})\ge \delta/4$ and if
$a=\{n_1,\ldots,n_k\}$ where $n_1<n_2<\cdots<n_k$ then we have
$$
|\{b\in a+:\psi(b)=m\}| =\cases 0 & \text{ if } m\le n_1+\cdots+n_k\\
1& \text{ if }m>n_1+\cdots+n_k.\\ \endcases
$$

Now for each $a\in T\setminus{\emptyset}$ we can choose $y_a^*\in
B_{Y^*}\cap E_{\psi(a)}^{\perp}$ so that $|y^*(y_a)|\ge \delta/4.$  Note
that the set $\{y_a^*:a\in T\}$ forms a weak$^*$-null sequence.  For
convenience let $y^*_{\emptyset}=0.$

Consider the closed unit interval $I=[0,1]$ and let $D$ be
the set of dyadic rationals $k/2^n$ where $1\le k\le 2^{n}-1$ and
$n\in\mathbb N.$ Let
$Z$ be the space of all real-valued functions
$f$ on
$I$ which are continuous on $I\setminus D$ and such that on $D$ both
left- and right-limits $f(q-)$ and $f(q+)$ exist
$$ f(q) =\frac12(f(q-)+f(q+)).$$
It is easy to see that $Z$ equipped with the sup-norm is isometric to
$C(\Delta)$ where $\Delta$ is the Cantor set.  Then $C(I)$ is a closed
subspace of $Z$ and $Z/C(I)\approx c_0(D)$ with the quotient map being
given by $Qf= (f(q+)-f(q-))_{q\in D}.$

Now we can define a one-one map $\varphi:T\to D$ with the property that
$\lim_{b\in a+}\varphi(b)=\varphi(a)$ and
$|\{b:\varphi(b)>\varphi(a)\}|=|\{b:\varphi(b)<\varphi(a)\}|=\infty$ for
$a\in T.$

Next define an operator $L:Y\to c_0(D)$ by putting $Ly(q)=y_a^*(y)$ if
$\varphi(a)=q$ and $Ly(q)=0$ otherwise; then $\|L\|\le 1.$  Then
$LQ_Y:\ell_1\to c_0(D)$ can be lifted to an operator $U:\ell_1\to Z$ so
that $QU=LQ_Y$ and $\|U\|\le 2.$   Then $U$ maps $X$ into $C(I)$ and by
assumption this restriction $U|_X$ has an extension $V:\ell_1\to C(I)$
with $\|V\|\le 2\lambda.$ Now $U-V$ factors to an operator $U-V=RQ_Y$
where $R:Y\to Z$ satisfies $\|R\|\le 2(\lambda+1)$ and $QR=L.$

  We can then write $R$ in the form
$$ Ry(q)= \langle y,h(q)\rangle$$ where $h:I\to Y^*$ is
weak$^*$-continuous except on points of $D$ and has left- and right-
weak$^*$-limits $h(q-)$ and $h(q+)$ on $D$ with
$$ h(q) =\frac12(h(q-)+h(q+)).$$
Note that $\|h(q)\|\le \lambda +1$
$$ h(q+)-h(q-) = \cases y_a^* &\text{ if } q=\varphi(a)\\ 0 &\text{ if }
q\notin \varphi(T).\endcases $$

Finally we build a branch $\beta=\{a_0,a_1,\ldots\}$ so that for each
$n\ge 1$ there exists $y_n^*=h(\varphi(a_n)+)$ or
$y_n^*=h(\varphi(a_n)-)$ so that
$$|\langle y_{a_k},y_n^*\rangle|>\delta/10$$
 for $1\le k\le n.$
This is done by induction.  Let $a_0=\emptyset$ and $a_1$ be any element
of $T$ with $|a_1|=1.$  Then since
$$ \frac{\delta}{4} \le \langle y_{a_1},y^*_{a_1}\rangle =\langle
y_{a_1},h(\varphi(a_1)+)-h(\varphi(a_1)-)\rangle $$
we can choose an appropriate sign so that the inductive hypothesis holds
when $n=1.$
Now suppose $a_0,\ldots,a_{n-1}$ have been chosen and that
$$ |\langle y_{a_k},y_{n-1}^*\rangle|
>\frac{\delta}{10}$$
for $1\le k\le n-1.$  Let us assume that
$y_{n-1}^*=h(\varphi(a_{n-1})+)$; the other case is similar. Then there
exists
$\eta>0$ so that if
$\varphi(a_{n-1})<q<\varphi(a_{n-1})+\eta$ we have for some $\rho>0,$
$$ |\langle y_{a_k},h(q)\rangle| >\frac{\delta}{10}+\rho$$
for $1\le k\le n-1.$  Then we can choose $a_n\in a_{n-1}+$ so that
$\varphi(a_{n-1})<\varphi(a_n)<\varphi(a_n)+\eta.$
Then
$$ |\langle y_{a_k}, h(\varphi(a_n)\pm)\rangle|>\frac{\delta}{10}$$ for
$1\le k\le n-1.$
Now $$\langle y_{a_n},h(\varphi(a_n)+)-h(\varphi(a_n)-)\rangle \ge
\frac{\delta}{4}$$
so that we can choose $y_n^*=h(\varphi(a_n)\pm)$ to satisfy the inductive
hypothesis.
This completes the inductive construction of the branch $\beta.$  Finally
we let $y^*$ be any weak$^*$-cluster point of the sequence
$((2\lambda+2)^{-1}y_n^*)_{n=1}^{\infty}$ so that $\|y^*\|\le 1$ and
$$|y^*(y_a)|\ge  \frac{\delta}{20(\lambda+1)}$$ for all $a\in \beta.$
This shows that $Y$ has the very strong Schur property with constant
$1/20(\lambda+1).$
\end{proof}

Our next theorem is then a partial converse of the Johnson-Zippin theorem
of \cite{JZ}.

\begin{Thm}\label{JZ2} Suppose $X$ is a closed subspace of $X$ such that
$(X,\ell_1)$ has (EP) and one of the following holds:
\newline
(i) $\ell_1/X$ has a (UFDD).\newline
(ii) $\ell_1/X$ is isomorphic to the dual of a Banach space $Y$ which
embeds into a space with a (UFDD).
\newline
Then $\ell_1/X$ is isomorphic to the dual of a subspace of $c_0$ and
there
is automorphism $\tau$ of $\ell_1$ so that $\tau(X)$ is weak$^*$-closed.
\end{Thm}

\begin{proof}
First note that (i) implies (ii).   In fact by Theorem \ref{JZSchur}
$\ell_1/X$ is a Schur space and hence any (UFDD) is boundedly complete so
that $\ell_1/X$ is a dual of a space with (UFDD).  If we assume (ii) then
Theorem \ref{c0dual} and  Theorem \ref{JZSchur} together yield the
result.\end{proof}

\noindent{\it Remark.} We can replace (i) by the assumption that
$\ell_1/X$
has (UMAP) (in some equivalent norm).  Indeed if $\ell_1/X$ has (UMAP) it
is shown in \cite{GK} that it has {\it commuting} (UMAP) and hence by
\cite{N}, Lemma 5.2, $\ell_1/X$ is the dual of a space with (UMAP).
Hence by Lemma \ref{unctype} and the remarks after Theorem \ref{c0dual}
we obtain that $\ell_1/X$ is the dual of a subspace of $c_0.$  As
observed already the classical results of \cite{LR} yield the existence
of the desired automorphism.

Let us also remark that, in the case when $\ell_1/X$ has a (UFDD) one can
easily deduce (from, say, results of \cite{JZ2}) that $\ell_1/X$ is
isomorphic to an $\ell_1-$sum of finite-dimensional spaces.

Note that we have proved:

\begin{Thm}\label{also} Let $X$ be a separable Banach space with a
(UFDD).
If $X$ has the very strong Schur property then $X$ is isomorphic to the
dual of a subspace of $c_0.$\end{Thm}

\enddocument